\documentclass[11pt]{article}

\usepackage{amsfonts}
\usepackage{enumerate}

\usepackage{amssymb}

\makeatletter
\renewcommand{\@biblabel}[1]{#1.\hfill}
\makeatother
\usepackage{amssymb}
\begin{document}

\centerline{\bf On a characterization theorem for the group of  $p$-adic numbers}

\bigskip

\centerline{\bf Gennadiy  Feldman}

\bigskip

\centerline{B.Verkin Institute for Low
Temperature  Physics and Engineering}

\centerline{ of the National Academy of
Sciences of Ukraine, Kharkov, Ukraine}

\centerline{
 feldman@ilt.kharkov.ua}

\bigskip

\begin{abstract}

It is well known  Heyde's characterization of the Gaussian
distribution on the real line: Let
 $\xi_1, \xi_2,\dots, \xi_n$, $n\ge 2,$ be independent random variables, let
$\alpha_j, \beta_j$ be nonzero constants such that
$\beta_i\alpha_i^{-1} + \beta_j\alpha_j^{-1} \ne 0$ for all $i \ne
j$.  If the conditional distribution of the linear form $L_2 =
\beta_1\xi_1 + \beta_2\xi_2+ \cdots + \beta_n\xi_n$  given $L_1 =
\alpha_1\xi_1 + \alpha_2\xi_2+\cdots + \alpha_n\xi_n$ is symmetric,
then all random variables $\xi_j$ are Gaussian. We prove an analogue
of this theorem for two independent random variables in the case
when they take values in the group of $p$-adic numbers $\Omega_p$,
and coefficients of linear forms are topological automorphisms of
$\Omega_p$.

\end{abstract}


\bigskip


\section{Introduction}

It is well known  Heyde's characterization of the Gaussian
distribution on the real line
  by the symmetry of the conditional distribution of one linear form given
another (\cite{He}, see also \cite[\S\,13.4.1]{KaLiRa}):

\medskip

{\bf Theorem A}. {\it Let
 $\xi_1, \xi_2,\dots, \xi_n$, $n\ge 2,$ be independent random variables, let
$\alpha_j, \beta_j$ be nonzero constants such that
$\beta_i\alpha_i^{-1} + \beta_j\alpha_j^{-1} \ne 0$  for all $i
\ne j$. If the conditional distribution of the linear form $L_2 =
\beta_1\xi_1 + \beta_2\xi_2+ \cdots + \beta_n\xi_n$ given $L_1 =
\alpha_1\xi_1 + \alpha_2\xi_2+\cdots + \alpha_n\xi_n$ is symmetric,
then all random variables $\xi_j$ are Gaussian.}

\medskip

Let $X$ be a second countable locally compact Abelian group, ${\rm
Aut}(X)$ be the group of topological automorphisms of $X$, $\xi_j,$
$j = 1, 2, \dots, n$, $n \ge 2$, be independent random variables
with values in $X$ and   distributions $\mu_j$. Consider the linear
forms $L_1 = \alpha_1\xi_1 + \alpha_2\xi_2+\cdots + \alpha_n\xi_n$
and $L_2 = \beta_1\xi_1 + \beta_2\xi_2+ \cdots + \beta_n\xi_n$,
where $\alpha_j, \beta_j \in {\rm Aut}(X)$. We  formulate  the
general problem.

\medskip

{\bf Problem  1.} {\it Let $X$ be a given group. Describe
distributions $\mu_j$ assuming that the conditional  distribution of
$L_2$ given $L_1$ is symmetric. In particular, for which groups $X$
the symmetry of the conditional distribution of $L_2$ given $L_1$
implies that all $\mu_j$ are either Gaussian distributions or belong
to a class of distributions which we can consider as a natural
analogue of the class of Gaussian distributions.}

\medskip

This problem for different classes of locally compact Abelian groups
was studied in series articles (see \cite{Fe2}-- \cite{Fe4},
\cite{Fe20bb}, \cite{Fe6}, \cite{Mi1}--\cite{MiFe1}).  In this  article we continue
these investigations. We study Problem 1 for two independent random
variables with values in the group of $p$-adic numbers.

We will use some results of the duality theory for locally compact
Abelian groups (see e.g. \cite{HeRo}). Before we formulate the main
theorem remind some definitions and agree on notation. For an
arbitrary locally compact Abelian
 group $X$  let $Y=X^\ast$ be its character group,
and  $(x,y)$  be the value of a character $y \in Y$ at an element $x
\in X$. If $K$ is a closed subgroup of $X$, we denote by $A(Y, K) =
\{y \in Y: (x, y) = 1$ for all $x \in K \}$ its annihilator. If
$\delta : X \mapsto X$ is a continuous endomorphism,
  then the adjoint endomorphism $\widetilde \delta : Y \mapsto Y$
is defined by the formula $(x, \widetilde \delta y) = (\delta x, y)$
for all $x \in X$, $y \in Y$. We note that $\delta \in {\rm Aut}(X)$
if and only if $\widetilde\delta \in {\rm Aut}(Y)$. Denote by $I$
the identity automorphism of a group.

Let  ${M^1}(X)$ be the convolution semigroup of probability
distributions on $X.$ For a distribution $\mu \in {M^1}(X)$ denote
by $\widehat \mu(y) = \int_X (x, y) d \mu(x)$ its characteristic
function. If $H$ is a closed subgroup of $Y$ and $\widehat \mu(y) =
1$ for $y \in H$, then $\widehat\mu(y+h) = \widehat\mu(y)$ for all
$y \in Y$, $h \in H$. For $\mu \in {M^1}(X)$, we define the
distribution $\bar \mu \in M^1(X)$ by the formula $\bar \mu(E) =
\mu(-E)$ for any Borel set $E \subset X$. Observe that $\widehat
{\bar \mu}(y) = \overline{\widehat \mu(y)}$. Let $K$ be a compact
subgroup of  $X$. Note that the characteristic function of the Haar
distribution $m_K$ is of the form
 \begin{equation}\label{1a}
 \widehat m_K(y) =
\left\{%
\begin{array}{ll}
    1, & \hbox{$y \in A(Y,K)$;} \\
    0, & \hbox{$y \notin A(Y,K)$.} \\
\end{array}%
\right.
\end{equation}
Denote by $I(X)$ the set of the idempotent distributions on $X$,
i.e. the set of shifts of the Haar distributions $m_K$ of compact
subgroups $K$ of $X$.

\section{The main theorem}

Let $p$ be a prime number. We need some properties of the group of
$p$-adic numbers $\Omega_p$ (see e.g. \cite[\S 10]{HeRo}). As a set
$\Omega_p$ coincides with the set of sequences of integers of the
form
 $x=(\dots,x_{-n}, x_{-n+1},\dots,
x_{-1}, x_0, x_1,\dots,x_n,\dots),$ where $x_n \in\{0, 1,\dots,
p-1\}$, such that  $x_n=0$ for $n < n_0$, where the number $n_0$
depends on $x$. Correspond to each element $x \in \Omega_p$  the
series
 $\sum\limits_{k=-\infty}^{\infty} x_k p^k.$ Addition and  multiplication of
 series are defined in a natural  way and define
the operations of addition and  multiplication in  $\Omega_p$. With
respect to these operations $\Omega_p$ is a field. Denote by
$\Delta_p$ a subgroup of $\Omega_p$ consisting of  $x \in \Omega_p$
such that $x_n=0$ for $n < 0$. This subgroup is called the group of
$p$-adic integers. Elements of the group $\Delta_p$ we   write in
the form $x=(x_0, x_1,\dots,x_n,\dots)$. The family of subgroups
 $\{p^m\Delta_p\}_{m=-\infty}^{\infty}$ can be considered as an open basis at zero of
 the group $\Omega_p$ and defines a topology on $\Omega_p$.
With respect to this topology the group $\Omega_p$ is locally
compact, noncompact, and totally disconnected. The character group
$\Omega_p^{*}$ of the group $\Omega_p$ is topologically isomorphic
to $\Omega_p$, and the value of a character $y \in \Omega_p^{*}$ at
an element $x \in \Omega_p$ is defined by the formula
 $$
(x, y)= \exp\Big[2\pi i\Big(\sum_{n=-\infty}^\infty
x_n\Big(\sum_{s=n}^\infty y_{-s}p^{-s+n-1}\Big)\Big)\Big]
$$
(\cite[(25.1)]{HeRo}).

Each automorphism $\alpha \in {\rm Aut}(\Omega_p)$ is the
multiplication by an element $x_\alpha \in \Omega_p,$ $x_\alpha \ne
0$, i.e. $\alpha g=x_{\alpha} g,$ $g \in \Omega_p.$
 If $\alpha \in
{\rm Aut}(\Omega_p),$ in order not to complicate notation,
 we will identify the automorphism
 $\alpha$ with the corresponding element
 $x_\alpha$, i.e. when we write $\alpha g$, we will suppose that  $\alpha
\in \Omega_p$. We note that $\widetilde\alpha=\alpha$. Denote by
$\Delta_p^0$ the subset of $\Delta_p$, consisting of all invertible
elements of $\Delta_p$, $\Delta_p^0=\{x=(x_0, x_1, \dots, x_n,\dots)
\in \Delta_p: x_0 \ne 0 \}$. Obviously, each element  $g \in
\Omega_p$ is represented in the form $g = p^k c,$ where $k$ is an
integer, and
 $c \in \Delta_p^0$. It is obvious that the multiplication by $c$  is a
topological automorphism of the group $p^m\Delta_p$ for any integer
$m$.

\medskip

For a fixed prime $p$ denote by ${\mathbf Z}(p^\infty)$ the set of
rational numbers of the form $${\left\{\frac{k}{p^n} : k=0, 1,
\dots,p^n-1, n=0,1,\dots\right\}}$$ and
 define the operation in ${\mathbf Z}(p^\infty)$ as
addition modulo 1. Then ${\mathbf Z}(p^\infty)$ is transformed into
an Abelian group, which we consider in the discrete topology.
Obviously,  this group is topologically isomorphic to the
multiplicative group of  $p^n$th roots of unity, where $n$ goes
through the nonnegative integers, considering in the discrete
topology. For a fixed $n$ denote by ${\mathbf Z}(p^n)$ the subgroup
of ${\mathbf Z}(p^\infty)$ consisting of all elements of the form
${\left\{\frac{k}{p^n} : k=0, 1, \dots,p^n-1\right\}}$. The group
 ${\mathbf Z}(p^n)$ is isomorphic to
the multiplicative group of  $p^n$th roots of unity.

\medskip

Let $\alpha \in {\rm Aut}(\Omega_p),$ $\alpha=p^kc,$ where $k \ge
0,$ $c \in \Delta_p^0$.  Let $l$ be an integer. It is easy to see
that
 $\alpha$ induces an epimorphism $\bar\alpha$ and $c$ induces an automorphism $\bar c$ on  the factor-group
$\Omega_p/p^l\Delta_p.$ Let $$x= (\dots,x_{-n}, x_{-n+1}, \dots,
x_{-1}, x_0, x_1,\dots,x_n,\dots) \in \Omega_p.$$ Define the mapping
$\tau:\Omega_p/p^l\Delta_p\mapsto {\mathbf{Z}}(p^\infty)$ by the
formula
 \begin{equation}\label{11a}
\tau(x+p^l\Delta_p)=\sum_{n=-\infty}^{l-1}x_np^{n-l}, \quad
x+p^l\Delta_p\in \Omega_p/p^l\Delta_p.
\end{equation}
Then $\tau$ is a topological isomorphism of the groups
$\Omega_p/p^l\Delta_p$ and $ {\mathbf{Z}}(p^\infty)$. Put
\begin{equation}\label{12a}
\widehat\alpha=\tau\bar\alpha\tau^{-1}, \quad \widehat c=\tau\bar
c\tau^{-1},
\end{equation}
and observe that $\widehat\alpha=p^k\widehat c$, $\widehat c\in {\rm
Aut}({\mathbf{Z}}(p^\infty))$, and $\widehat\alpha$ is an
epimorphism. If $c=(c_0, c_1, \dots, c_n,\dots) \in \Delta_p^0$,
then the automorphism $\widehat c$ acts in the following way. Put
$s_n=c_0+c_1p+c_2p^2+\dots+c_{n-1}p^{n-1}$. The restriction of the
automorphism $\widehat c$ to the subgroup ${\mathbf{Z}}(p^n) \subset
{\mathbf{Z}}(p^\infty)$ is of the form $\widehat cy = s_n y,$ $y \in
{\mathbf{Z}}(p^n)$, i.e. $\widehat c$ acts in ${\mathbf{Z}}(p^n)$ as
the multiplication by $s_n$.

We note that since $\Omega_p$ is a totally disconnected group, the
Gaussian distributions on   $\Omega_p$ are degenerated (\cite[Ch.
4]{P}), and the class of idempotent distributions on $\Omega_p$ can
be considered as a natural analog of the class of Gaussian
distributions.

\medskip

We    consider the case of two independent random variables
 $\xi_1$  and $\xi_2$ with values in the group $X=\Omega_p$,
i.e. we   consider the linear forms
 $L_1=\alpha_1\xi_1 + \alpha_2\xi_2$ and $L_2=\beta_2\xi_1 +
\beta_2\xi_2,$ where $\alpha_j, \beta_j \in {\rm Aut}(X)$. Assume
that the conditional distribution of the linear form
$L_2=\beta_2\xi_1 + \beta_2\xi_2$ given $L_1=\alpha_1\xi_1 +
\alpha_2\xi_2$ is symmetric. Our aim is to describe the possible
distributions  $\mu_j$ depending on $\alpha_j, \beta_j$. Introducing
the new independent random variables $\xi'_j=\alpha_j\xi_j$, $j=1,
2$, we can suppose  that $L_1=\xi_1 + \xi_2$ and $L_2=\delta_1\xi_1
+ \delta_2\xi_2,$ where $\delta_j \in {\rm Aut}(X).$ Note also that
the conditional distribution of the linear form $L_2=\delta_1\xi_1 +
\delta_2\xi_2$ given $L_1=\xi_1 + \xi_2$ is symmetric if and only if
the conditional distribution of $\delta L_2$, where $\delta\in {\rm
Aut}(X)$,  given $L_1$ is symmetric. It follows from this that
  we may assume that $L_2=\xi_1 + \alpha\xi_2$, where $\alpha \in {\rm Aut}(X)$.

We can formulate now the main result of the article.

\medskip

{\bf Theorem 1}. {\it Let $X=\Omega_p$. Let $\alpha=p^kc$, where $k\in \mathbf{Z}$ and $c=(c_0, c_1,\dots, c_n,\dots)
 \in \Delta_p^0$, be an arbitrary topological automorphism of the group $X$. Then the following statements hold.

{$1$.} Let $\xi_1$ and $\xi_2$ be independent random variables with
values in $X$ and distributions $\mu_1$ and $\mu_2$. Assume that the
conditional distribution of the linear form $L_2=\xi_1 + \alpha
\xi_2$ given $L_1=\xi_1+\xi_2$ is symmetric.  Then

$1(i)$ If  $p>2$, $k=0$ and $c_0 \ne p-1$, then $\mu_j=m_K*E_{x_j}$,
where $K$ is a compact subgroup of
 $X$ and $x_j\in X$. Moreover, if $c_0=1$, then $\mu_1$ and $\mu_2$
are degenerate distributions.

$1(ii)$ If $p=2$, $k=0$, $c_0=1$ and $c_1=0$, then $\mu_1$ and
$\mu_2$ are degenerate distributions.

$1(iii)$ If  $p>2$ and $|k|=1$, then either  $\mu_1\in I(X)$ or
$\mu_2\in I(X)$.

{$2$.} If one of the following conditions holds:

$2(i)$ $p>2,$ $k=0,$ $c_0=p-1;$

$2(ii)$ $p=2,$ $k=0,$ $c_0=c_1=1$;

$2(iii)$ $p=2,$ $|k|=1$;

$2(iv)$ $p\ge 2,$ $|k| \ge 2$,

\noindent then there exist independent random variables $\xi_1$ and
$\xi_2$  with values in $X$ and distributions $\mu_1$ and $\mu_2$
such that the conditional distribution of the linear form $L_2=\xi_1
+ \alpha \xi_2$ given $L_1=\xi_1+\xi_2$ is symmetric whereas $\mu_1,
\mu_2 \notin I(X)$.}

\medskip

It is easy to see that cases $1(i)$--$1(iii)$ and  $2(i)$--$2(iv)$
exhaust all possibilities for $p$ and $\alpha$.

\section{Lemmas}

To prove Theorem 1 we need some lemmas.

\medskip

{\bf Lemma 1} (\cite{Fe4}, see also \cite[\S 16]{Fe5}). {\it Let $X$
be a second countable locally compact Abelian group, $Y$ be its
character group. Let $\xi_1$ and $\xi_2$ be independent random
variables with  values in $X$ and distributions $\mu_1$ and $\mu_2$.
Let $\alpha_j, \beta_j$ be continuous endomorphisms of $X$. The
conditional distribution of
 the linear form $L_2=\beta_1\xi_1 + \beta_2\xi_2$ given  $L_1=\alpha_1\xi_1 + \alpha_2\xi_2$
is symmetric if and only if the characteristic functions
$\widehat\mu_j(y)$ satisfy   the equation}
 \begin{equation}\label{2}
\widehat \mu_1(\widetilde\alpha_1u + \widetilde\beta_1  v)\widehat
\mu_2(\widetilde\alpha_2u + \widetilde\beta_2  v) = \widehat
\mu_1(\widetilde\alpha_1u - \widetilde\beta_1  v)\widehat
\mu_2(\widetilde\alpha_2u - \widetilde\beta_2  v), \quad u, v \in Y.
\end{equation}

\medskip

Observe that in \cite{Fe4} this lemma is proved when $\alpha_j,
\beta_j\in {\rm Aut}(X)$, but the proof is valid without any changes
for the case when $\alpha_j, \beta_j$ are continuous endomorphisms.

\medskip

{\bf Lemma 2}. {\it Let $X = \Omega_p$,
  $\xi_1$ and $\xi_2$ be independent random variables with  values in $X$ and
  distributions $\mu_1$ and $\mu_2$ such that
$\mu_j(y) \ge 0$, $j = 1, 2$. Let $\alpha=p^kc \in {\rm Aut}(X)$,
where $k\in \mathbf{Z}$ and  $c \in \Delta_p^0$. Assume also that
$\alpha\ne -I$. If the conditional distribution of the linear form
$L_2=\xi_1 + \alpha \xi_2$ given $L_1=\xi_1+\xi_2$ is symmetric,
then there exists a closed subgroup $H\subset Y$, such that
$\widehat\mu_j(y)=1$ for $y \in H$, $j=1, 2$.}

\medskip

{\bf Proof}. Taking into account  that the conditional distribution
of the linear form $L_2$ given $L_1$ is symmetric if and only if the
conditional distribution of the linear form $\alpha^{-1} L_2$
  given $L_1$ is symmetric,
  we may assume without loss of generality that $k\ge 0$.
Since $X=\Omega_p$, we have $\widetilde\alpha=\alpha$, and hence
equation (\ref{2}) takes the form
\begin{equation}\label{3n}
\widehat \mu_1(u+v)\widehat \mu_2(u+\alpha v)=\widehat
\mu_1(u-v)\widehat \mu_2(u-\alpha v), \quad u, v \in Y.
\end{equation}
If $\alpha = I$, then put in equation (\ref{3n}) $u=v=y$. We get
that $H=Y$, i.e. $\mu_1,$ and $\mu_2$ are degenerate distributions.
So, we will assume that  $\alpha \ne I$.
 Since $\widehat\mu_1(0)=\widehat\mu_2(0)=1$, we can choose a neighborhood
at zero $V$ of $Y$  such that $\widehat\mu_j(y) >0$ for $y \in V$,
$j=1, 2$. Obviously, we can assume that   $V=p^l\Delta_p$ for some
$l$. Since $k\ge 0$, we have $\alpha(p^l\Delta_p)\subset
p^l\Delta_p$.   Put $\psi_j(y) = - \log \widehat\mu_j(y)$, $y \in
p^l\Delta_p$, $j=1, 2$. It follows from (\ref{3n}) that
 the functions $\psi_1(y)$ and $\psi_2(y)$
satisfy the equation
\begin{equation}\label{4}
\psi_1(u + v)+\psi_2(u + \alpha v)  - \psi_1(u - v)-\psi_2(u -
\alpha v)=0, \quad u, v \in p^l\Delta_p.
\end{equation}
We use the finite difference method to solve equation (\ref{4}). Let
$\psi(y)$ be a function on  $Y$, and $h$ be an arbitrary element of
 $Y$. Denote by $\Delta_h$
the finite difference operator
$$\Delta_h \psi(y)=\psi(y+h)-\psi(y).$$
Let $k_1$ be an arbitrary element of $p^l\Delta_p$. Put $h_1 =\alpha
k_1$ and hence, $h_1 -\alpha k_1 = 0$. Substitute $u+h_1$ for $u$
and $v+k_1$ for $v$  in equation (\ref{4}) Subtracting equation
(\ref{4}) from the obtained equation we find
\begin{equation}\label{5}
\Delta_{l_{11}}\psi_1(u+v) + \Delta_{l_{12}}\psi_2(u+\alpha v) -
\Delta_{l_{13}}\psi_1(u-v)=0, \quad u, v \in p^l\Delta_p,
\end{equation}
\noindent where $l_{11}= (\alpha + I)k_1,$ $ l_{12}=2 \alpha k_1,$ $
l_{13}= (\alpha-I)k_1.$ Let $k_2$ be an arbitrary element of
$p^l\Delta_p$. Put $h_2 =k_2$ and hence, $h_2 -k_2 = 0$. Substitute
$u+h_2$ for $u$ and $v+k_2$ for $v$  in equation (\ref{5}).
 Subtracting equation (\ref{5}) from the obtained
equation we arrive at
\begin{equation}\label{6}
\Delta_{l_{21}}\Delta_{l_{11}}\psi_1(u+v) +
\Delta_{l_{22}}\Delta_{l_{12}}\psi_2(u+\alpha v) = 0,
 \quad u, v \in p^l\Delta_p,
\end{equation}
where $l_{21}=2 k_2,$ $l_{22}= (I+\alpha)k_2.$ Let $k_3$ be an
arbitrary element of $p^l\Delta_p$. Put $h_3 =-\alpha k_3$ and
hence, $h_3 +\alpha k_3 = 0$. Substitute $u+h_3$ for $u$ and $v+k_3$
for $v$  in equation (\ref{6}).
 Subtracting equation (\ref{6}) from the obtained
equation
 we find
\begin{equation}\label{7}
\Delta_{l_{31}}\Delta_{l_{21}}\Delta_{l_{11}}\psi_1(u+v) = 0,
 \quad u, v \in p^l\Delta_p,
\end{equation}
where $l_{31}= (I-\alpha)k_3.$ Substituting $v=0$ into (\ref{7}),
taking into account the expressions for  $l_{11}$, $l_{21}$,
$l_{31}$ and the fact that $k_1, k_2, k_3$   are arbitrary elements
of $p^l\Delta_p$, we find from (\ref{7}) that there exists a
subgroup $p^m\Delta_p\subset p^l\Delta_p$, where the function
$\psi_1(y)$ satisfies the equation
\begin{equation}\label{8}
\Delta_h^3\psi_1(y)=0, \quad h, y \in p^m\Delta_p.
\end{equation}
 Taking into account that
$p^m\Delta_p$ is a compact group, we conclude from (\ref{8}) that
$\psi_1(y)=const$, $y \in p^m\Delta_p$ (\cite[\S 5]{Fe5}). Since
$\psi_1(0)=0$, we have $\psi_1(y)= 0$, for $y \in p^m\Delta_p$. This
implies that $\widehat\mu_1(y)=1$ for $y \in p^m\Delta_p$. For the
distribution $\mu_2$ we reason similarly and find a subgroup
$p^n\Delta_p$ such that $\widehat\mu_2(y)=1$ for $y\in p^n\Delta_p$.
Put $H=p^m \Delta_p \cap p^n \Delta_p$. Lemma 2 is proved.

\medskip

{\bf Lemma 3} (\cite{Fe2}, see also \cite[\S 16]{Fe5}). {{\it  Let
$X$ be a finite Abelian group containing no elements of order $2$,
 $\alpha$ be an automorphism of
$X$ such that  $I\pm \alpha\in {\rm Aut}(X)$. Let $\xi_1$ and
$\xi_2$ be independent random variables with values in
 the group $X$ and  distributions  $\mu_1$ and $\mu_2$.
If the conditional distribution of the linear form $L_2=\xi_1 +
\alpha \xi_2$ given $L_1 = \xi_1 + \xi_2$ is symmetric, then
$\mu_j=m_K*E_{x_j}$, where $K$ is a finite subgroup of
 $X$  and $x_j\in X$, $j=1, 2$}.

\section{Proof of statements $1(i)-1(iii)$}

Before we begin the proof make the following remarks. It is obvious
that the characteristic functions of the distributions $\bar\mu_j$
also satisfy equation (\ref{3n}). This implies that the
characteristic functions of the distributions
 $\nu_j=\mu_j*\bar\mu_j$ satisfy equation (\ref{3n}) too.  We have
$\widehat\nu_j(y)=|\widehat\mu_j(y)|^2 \ge 0$, $j = 1, 2$. Hence,
when we prove statements  { $1(i)$}--{ $1(iii)$} we may assume
without loss of generality that
 $\mu_j(y)
\ge 0$, $j = 1, 2$, because $\mu_j$ and $\nu_j$ are  idempotent
distributions or degenerate distributions simultaneously.

Taking into account  that the conditional distribution of the linear
form $L_2$ given $L_1$ is symmetric if and only if the conditional
distribution of the linear form $\alpha^{-1} L_2$
  given $L_1$ is symmetric,
  we may assume without loss of generality that $k\ge 0$.

Put $f(y)=\widehat\mu_1(y)$, $g(y)=\widehat\mu_2(y)$ and write
equation (\ref{3n}) in the form
\begin{equation}\label{3}
f(u+v)g(u+\alpha v)=f(u-v)g(u-\alpha v), \quad u, v \in Y,
\end{equation}
where  $\alpha \in {\rm Aut}(Y)$, $\alpha= p^k c$, $k \ge 0$, $c \in
\Delta_p^0$. We also assume that $f(y) \ge 0$, $g(y) \ge 0$. In fact
we will study solutions of equation  (\ref{3}).

\medskip

{\bf Proof of statements $1(i)$ and $1(ii)$.} By Lemma 1, the
characteristic functions $f(y)$ and $g(y)$ satisfy equation
(\ref{3}). Put
\begin{equation}\label{9}
E = \{ y \in Y: f(y) = g(y) = 1 \}.
\end{equation}
If either $\mu_1$ or $\mu_1$ is a nondegenerate distribution, then
$E \ne \Omega_p$. Observe that in cases $1(i)$ and $1(ii)$
$\alpha\ne -I$. Then by Lemma 2, $E \ne\{0\}.$ Thus, $E$ is a
nonzero proper closed subgroup of
 $\Omega_p,$ and hence $E=p^l\Delta_p$ for some $l$.
It follows from (\ref{9})  that
\begin{equation}\label{13a}
f(y + h) = f(y), \quad g(y + h) = g(y), \quad y \in Y, \ h \in
p^l\Delta_p.
\end{equation}
Taking into account (\ref{13a}), denote by $\bar f(y)$ and $\bar
g(y)$ the functions induced by  the functions $f(y)$ and $g(y)$
 on the factor-group $Y/p^l\Delta_p$.   Put $\widehat f=\bar f\circ \tau^{-1}$,
where $\tau$ is defined by formula (\ref{11a}). We get from
(\ref{3}) that the functions $\widehat f(y)$ and $\widehat g(y)$
satisfy the equation
\begin{equation}\label{10}
\widehat f(u+v)\widehat  g(u+\widehat\alpha v)=\widehat
f(u-v)\widehat  g(u-\widehat\alpha v), \quad u, v \in {\mathbf
Z}(p^\infty),
\end{equation}
where $\widehat\alpha$ is defined by formula (\ref{12a}).
 It follows from (\ref{9}) that
\begin{equation}\label{11}
\{ y \in {\mathbf Z}(p^\infty): \widehat f(y) = \widehat g(y) = 1 \}
= \{0\}.
\end{equation}

\medskip

{\bf Statement $1(i)$}.  Assume first that $c_0 \ne 1$. Since
$p>2$, $k=0$ and  $c_0 \ne p-1$, we have
\begin{equation}\label{13}
I \pm \widehat\alpha \in {\rm Aut}({\mathbf{Z}}(p^\infty)).
\end{equation}
We note that for any $n$ the restriction of any automorphism of the
group ${\mathbf{Z}}(p^\infty)$ to the subgroup ${\mathbf{Z}}(p^n)
\subset {\mathbf{Z}}(p^\infty)$ is an automorphism of the subgroup
${\mathbf{Z}}(p^n).$  Consider the restriction of   equation
(\ref{10}) to the subgroup
 ${\mathbf{Z}}(p^n)$. Observe that
$({\mathbf{Z}}(p^n))^*\cong{\mathbf{Z}}(p^n)$, and the group
${\mathbf{Z}}(p^n)$ contains no elements of order 2. Taking into
account that (\ref{13}) holds, we can apply Lemmas 1 and  3 to the
group ${\mathbf{Z}}(p^n)$ and get that the restrictions of the
characteristic functions $\widehat f(y)$ and $\widehat g(y)$ to the
subgroup ${\mathbf{Z}}(p^n)$ take only two values $0$ and $1$.
Moreover, $\widehat f(y)=\widehat g(y)$ for $y\in
{\mathbf{Z}}(p^n)$. Hence, the characteristic functions $\widehat
f(y)$ and $\widehat g(y)$ on the group ${\mathbf{Z}}(p^\infty)$ take
also only two values $0$ and $1$, and $\widehat f(y)=\widehat g(y)$
for $y\in {\mathbf{Z}}(p^\infty)$.
 Then the standard reasoning show that
 $\mu_j=m_K*E_{x_j}$, where $K$
is a compact subgroup of
 $X$ and $x_j\in X$.

\medskip

Assume now that  $c_0=1$. Since $k=0$, the restriction of the
automorphism $\widehat\alpha$ to ${\mathbf{Z}}(p) \subset
{\mathbf{Z}}(p^\infty)$ is the identity automorphism. Hence, the
restriction of equation (\ref{10}) to ${\mathbf{Z}}(p)$ takes the
form
\begin{equation}\label{12}
\widehat f(u+v) \widehat g(u+v)=\widehat f(u-v)\widehat g(u-v),
\quad u, v \in {\mathbf{Z}}(p).
\end{equation}
Substituting here $u=v=y$, we get $\widehat f(2y) \widehat g(2y)=1$,
$y \in {\mathbf{Z}}(p)$. Since $p>2$, this implies that $$\widehat
f(y)=\widehat g(y)= 1, \quad y \in {\mathbf{Z}}(p),$$ but this
contradicts to  (\ref{11}). Thus, $\mu_1$ and $\mu_2$ are degenerate
distributions.

\medskip

{\bf Statement $1(ii)$}.  Since $k=0$, $c_0=1$ and $c_1=0$,  the
restriction of the automorphism $\widehat\alpha$ to the subgroup
${\mathbf{Z}}(4) \subset {\mathbf{Z}}(2^\infty)$ is the identity
automorphism.  Hence, the restriction of equation (\ref{10}) to the
subgroup ${\mathbf{Z}}(4)$ takes the form (\ref{12}), where $u, v
\in {\mathbf{Z}}(4)$. Substituting in (\ref{12}) $u=v=y$, we get
$$\widehat f(2y) \widehat g(2y)=1, \quad y \in  {\mathbf{Z}}(4).$$
This implies that $\widehat f(y)=\widehat g(y)= 1$ for $y \in
{\mathbf{Z}}(2)$, but this contradicts to  (\ref{11}). Thus, $\mu_1$
and $\mu_2$ are degenerate distributions.

\medskip

The proof of statement $1(iii)$ is based on the proof of the
following proposition which is of interest in its own right.

\medskip

{\bf Proposition 1}. {\it Let $X = \Omega_p$,
  $\xi_1$ and $\xi_2$ be independent random variables with  values in $X$ and
  distributions $\mu_1$ and $\mu_2$.  Let $\alpha \in {\rm
Aut}(X)$, and $\alpha\ne -I$. If the conditional distribution of the
linear form  $L_2=\xi_1 + \alpha \xi_2$ given $L_1=\xi_1+\xi_2$ is
symmetric, then either $\mu_1$ and $\mu_2$ are degenerate
distributions or
 there exists a closed subgroup $M\subset Y$, such that $\widehat\mu_j(y)=0$ for $y \notin M, \ j=1, 2$.}

\medskip

{\bf Proof}. Let $\alpha=p^kc \in {\rm Aut}(X)$, where $k\in
\mathbf{Z}$ and  $c \in \Delta_p^0$. We can assume without loss of
generality that $k\ge 0$. Otherwise we consider the  linear form
$L_2=\alpha^{-1}\xi_1+\xi_2$ instead of $L_2=\xi_1 + \alpha \xi_2$.
Taking into account that in cases $1(i)$ and $1(ii)$ of Theorem 1
Proposition 1
 holds, it remains to consider
the following cases:

1. $p\ge 2$, $ k \ge 1$;

2. $p>2,$ $k=0,$ $c_0=p-1;$

3. $p=2,$ $k=0,$ $c_0=c_1=1$.

Obviously, we can assume without loss of generality that
 $\widehat\mu_j(y) \ge 0,$ $j=1, 2$. Put $f(y)=\widehat\mu_1(y)$, $g(y)=\widehat\mu_2(y)$.
We have $f(-y)=f(y)$ and $g(-y)=g(y)$. Reasoning as in the proof of
statements $1(i)$--$1(ii)$ and retaining the same notation we arrive
at equation (\ref{10}). Set $\beta=I-\alpha,$ $\gamma=I+\alpha$.
Substituting into equation (\ref{10}) $u=v=y$, we obtain
\begin{equation}\label{14}
\widehat f(2y)\widehat  g(\widehat\gamma y)=\widehat g(\widehat\beta
y), \quad y \in {\mathbf{Z}}(p^\infty).
\end{equation}
Substituting into equation (\ref{10})  $u=\widehat\alpha y,$ $v=y$,
we get
\begin{equation}\label{15}
\widehat f(\widehat\gamma y)\widehat g(2\widehat\alpha y) =\widehat
f(\widehat\beta y), \quad y \in {\mathbf{Z}}(p^\infty).
\end{equation}

\medskip

1. $p\ge 2$, $ k \ge 1$. Since $k \ge 1$, we have $\widehat\beta,
\widehat\gamma \in {\rm Aut}({\mathbf{Z}}(p^\infty))$. Put $\kappa=
\beta \gamma^{-1}$. Equations (\ref{14})  and (\ref{15})  imply
\begin{equation}\label{16}
\widehat g(\widehat\kappa y)= \widehat
f(2\widehat\gamma^{-1}y)\widehat g(y), \quad y \in
{\mathbf{Z}}(p^\infty),
\end{equation}
and
\begin{equation}\label{17}
\widehat f(\widehat\kappa y)=\widehat f(y) \widehat
g(2\widehat\alpha\widehat\gamma^{-1}y), \quad y \in
{\mathbf{Z}}(p^\infty).
\end{equation}
Since  $0 \le \widehat f(y) \le 1$, it follows from (\ref{16}) that
$$\widehat g(\widehat\kappa y)\le \widehat g(y), \quad y \in
{\mathbf{Z}}(p^\infty).$$ This implies that for any natural $n$ the
inequalities
\begin{equation}\label{18}
\widehat g(\widehat\kappa^ny) \le\dots\le \widehat g(\widehat\kappa
y) \le\widehat g(y), \quad y \in {\mathbf{Z}}(p^\infty),
\end{equation}
hold. Let $y \in {\mathbf{Z}}(p^\infty)$. Then $y \in
{\mathbf{Z}}(p^l)$ for some $l$, and hence $\widehat\kappa y \in
{\mathbf{Z}}(p^l).$ It follows from this that $\widehat\kappa^my=y$
for some $m$, depending generally speaking on $y$. Substituting in
(\ref{18}) $n=m$ we get
\begin{equation}\label{19}
\widehat g(\widehat\kappa y) =\widehat g(y), \quad  y\in
{\mathbf{Z}}(p^\infty).
\end{equation}
Reasoning similarly, we obtain from (\ref{17}) that
\begin{equation}\label{20}
\widehat f(\widehat\kappa y) =\widehat f(y), \quad  y\in
{\mathbf{Z}}(p^\infty).
\end{equation}

Assume that there exists a sequence of elements $y_n \in
{\mathbf{Z}}(p^\infty)$ such that:

$(a)$ the order of element $y_n$ is equal to $p^{i_n},$ $i_n
\rightarrow \infty$;

$(b)$ $\widehat g(y_n) \ne 0$.

Then it follows from (\ref{16}) and (\ref{19}) that,
\begin{equation}\label{21}
\widehat f(2\widehat\gamma^{-1}y_n) \ =1.
\end{equation}
Suppose that $p>2$. Then  (\ref{21})  implies that  $\widehat f(y) =
1$ for $y \in {\mathbf{Z}}(p^{i_n})$, because the order of the
element $2\widehat\gamma^{-1}y_n$ is equal to $p^{i_n}$, and hence
the element $2\widehat\gamma^{-1}y_n$ generates the subgroup
${\mathbf{Z}}(p^{i_n})$. If  $p=2$, then $\widehat f(y) = 1$ for $y
\in {\mathbf{Z}}(2^{i_n-1})$, because the order of the element
$2\widehat\gamma^{-1}y_n$ is equal to $2^{i_n-1}$, and the element
 $2\widehat\gamma^{-1}y_n$ generates the subgroup
${\mathbf{Z}}(2^{i_n-1})$. Thus, for $p\ge 2$ we have $\widehat f(y)
= 1$ for $ y \in {\mathbf{Z}}(p^\infty),$ and equation   (\ref{10})
implies that
$$
\widehat g(u+\widehat\alpha v)= \widehat g(u-\widehat\alpha v),
\quad u,v \in {\mathbf{Z}}(p^\infty).
$$
It follows from this that $\widehat g(2\widehat\alpha y)= 1$ for $y
\in {\mathbf{Z}}(p^\infty)$, and hence $\widehat g(y) = 1$ for $y
\in {\mathbf{Z}}(p^\infty),$ because $2\widehat\alpha$ is an
epimorphism. We proved that $\mu_1,$ and $\mu_2$ are degenerate
distributions.

The similar reasoning show that if there exists a sequence of
elements $z_n \in {\mathbf{Z}}(p^\infty)$ such that:

$(a)$ the order of element $z_n$ is equal to $p^{j_n},$ $j_n
\rightarrow \infty$;

$(b)$ $\widehat f(z_n) \ne 0$,

\noindent then $\mu_1,$ and $\mu_2$ are also degenerate
distributions.

From what has been said it follows that if $\mu_1,$ and $\mu_2$ are
nondegenerate distributions, then there exists $n$ such that
$\widehat f(y)=\widehat g(y)=0$ for
 $y \notin {\mathbf{Z}}(p^n)$. Proposition 1 in case 1
 follows directly from this.

\medskip

2. $ p>2,$ $k=0,$ $c_0=p-1.$ Since $p >2$, we have $\widehat\beta\in
{\rm Aut}({\mathbf{Z}}(p^\infty))$. We find from equation (\ref{14})
that
\begin{equation}\label{22}
\widehat g(y)=\widehat f(2 \widehat\beta^{-1}y) \widehat
g(\widehat\gamma \widehat\beta^{-1}y), \quad y \in
{\mathbf{Z}}(p^\infty).
\end{equation}
Since $0 \le \widehat g(y) \le 1$, we find from (\ref{22}) that
\begin{equation}\label{23}
\widehat g(y) \le \widehat f(2 \widehat\beta^{-1}y), \quad y \in
{\mathbf{Z}}(p^\infty).
\end{equation}
Reasoning similarly we get from  (\ref{15}) that
\begin{equation}\label{23a}
\widehat f(y) = \widehat
f(\widehat\gamma\widehat\beta^{-1}y)g(2\widehat\alpha
\widehat\beta^{-1}y), \quad y \in {\mathbf{Z}}(p^\infty).
\end{equation}
Taking into account that $0\le \widehat f(y)\le 1$, this implies
that
\begin{equation}\label{24}
\widehat f(y) \le \widehat g(2\widehat\alpha \widehat\beta^{-1}y),
\quad y \in {\mathbf{Z}}(p^\infty).
\end{equation}
Inequalities   (\ref{23})  and  (\ref{24})  imply the inequalities
\begin{equation}\label{25}
\widehat  g(y) \le \widehat  f(2 \widehat\beta^{-1}y) \le \widehat
g(4\widehat \alpha \widehat \beta^{-2}y), \quad \widehat  f(y) \le
\widehat  g(2\widehat \alpha \widehat\beta^{-1}y) \le \widehat
f(4\widehat\alpha\widehat\beta^{-2}y), \quad y \in
{\mathbf{Z}}(p^\infty).
\end{equation}
Reasoning as in the proof of case 1, we find from  (\ref{25})  that
\begin{equation}\label{14a}
\widehat g(y)= \widehat g(4\widehat\alpha \widehat\beta^{-2}y),
\quad \widehat f(y)=\widehat f(4\widehat\alpha \widehat\beta^{-2}y),
\quad y \in {\mathbf{Z}}(p^\infty).
\end{equation}
We find from  (\ref{25})  and  (\ref{14a})    that
\begin{equation}\label{26}
\widehat g(y) =\widehat f(2 \widehat\beta^{-1}y), \quad
 \widehat f(y) = \widehat g(2\widehat\alpha \widehat\beta^{-1}y), \quad y
\in {\mathbf{Z}}(p^\infty).
\end{equation}
It follows from   (\ref{26}) (\ref{23a}) and (\ref{22}) that if
$\widehat g(y_0) \ne 0$ for some  $y_0 \in {\mathbf{Z}}(p^\infty)$,
then $\widehat g(\widehat\gamma\widehat\beta^{-1}y_0)=1$, and if
$\widehat f(y_0) \ne 0$, then $\widehat
f(\widehat\gamma\widehat\beta^{-1}y_0)=1$. We complete the proof as
in case 1.

\medskip

3. $p=2,$ $k=0,$ $c_0=c_1=1$. Put $\beta=2\beta_1,$
$\gamma=2\gamma_1.$ Then $\widehat\beta_1 \in {\rm
Aut}({\mathbf{Z}}(2^\infty)),$ and $\widehat\gamma_1$ is an
epimorphism. It follows from (\ref{14})  that
\begin{equation}\label{27}
\widehat f(y)\widehat g(\widehat\gamma_1 y)=\widehat
g(\widehat\beta_1y), \quad y \in {\mathbf{Z}}(2^\infty).
\end{equation}
Similarly, we find from (\ref{15})  that
\begin{equation}\label{28}
\widehat f(\widehat\gamma_1y)\widehat g(\widehat\alpha y) =\widehat
f(\widehat\beta_1 y), \quad y \in {\mathbf{Z}}(2^\infty).
\end{equation}
It follows from  (\ref{27})  and (\ref{28}) that
$$
\widehat g(y)=\widehat f({\widehat\beta_1}^{-1} y)\widehat
g(\widehat\gamma_1{\widehat\beta_1}^{-1}y), \quad \widehat
f(y)=\widehat f(\widehat\gamma_1{\widehat\beta_1}^{-1} y)\widehat
g(\widehat\alpha{\widehat\beta_1}^{-1}y), \quad y \in
{\mathbf{Z}}(2^\infty).
$$
Hence,
$$
\widehat g(y) \le \widehat f({\widehat\beta_1}^{-1} y), \quad
\widehat f(y) \le \widehat g(\widehat\alpha{\widehat\beta_1}^{-1}
y).
$$
This implies that
 $$\widehat g(y) \le
\widehat f({\widehat\beta_1}^{-1} y) \le \widehat
g(\widehat\alpha{\beta_1}^{-2} y), \quad \widehat f(y) \le \widehat
g(\widehat\alpha{\widehat\beta_1}^{-1} y) \le \widehat
f(\widehat\alpha{\widehat\beta_1}^{-2} y).$$ We complete the proof
as in case 2. Proposition 1 is proved completely.

\medskip

{\bf Proof of statement $1(iii)$}. Observe that in case  $1(iii)$
 $\alpha\ne -I$. Reasoning as in the proof of statements $1(i)$ and $1(ii)$
and retaining the same notation we arrive at equation (\ref{10}).
Put
$$E_{\widehat f}=\{y \in {\mathbf{Z}}(p^\infty): \widehat f(y) \ne 0\}, \quad
 E_{\widehat g}=\{y \in {\mathbf{Z}}(p^\infty): \widehat g(y) \ne 0\},
$$
$$
B_{\widehat f}=\{y \in {\mathbf{Z}}(p^\infty): \widehat f(y)=1\},
\quad B_{\widehat g}=\{y \in {\mathbf{Z}}(p^\infty): \widehat g(y)
=1\}.
$$
Since $p>2$  and $k=1$, we have $2\widehat\gamma^{-1} \in {\rm
Aut}({\mathbf{Z}}(p^\infty))$. Moreover, since $k=1$, obviously,
(\ref{16}), (\ref{17}), (\ref{19}) and (\ref{20}) hold. Taking this
into account, we find from  (\ref{16}) and (\ref{19}) that
$E_{\widehat g} \subset B_{\widehat f}$. Analogously, we find from
(\ref{17}) and (\ref{20}) that $pE_{\widehat f} \subset B_{\widehat
g}$. If $\mu_2$ is a nondegenerate distribution, then $B_{\widehat
g}$ is a proper subgroup of ${\mathbf{Z}}(p^\infty)$, and hence
$B_{\widehat g}={\mathbf{Z}}(p^{n})$ for some  $n$. Since
$B_{\widehat g} \subset E_{\widehat g} \subset B_{\widehat f}$, it
follows from (\ref{11}) that
 $B_{\widehat  g} =\{0\}$. If $B_{\widehat  f} =\{0\}$, then
$E_{\widehat g} \subset B_{\widehat  f}  =\{0\}$, and hence $\mu_2$
is an idempotent distribution. If  $B_{\widehat  f} \ne \{0\}$, then
$pB_{\widehat f} \subset pE_{\widehat f} \subset B_{\widehat
g}=\{0\}$. This implies that $B_{\widehat f}= E_{\widehat f} =
{\mathbf{Z}}(p)$, and hence $\mu_1$ is an idempotent distribution.

\section{Proof of statements $2(i)$--$2(iv)$}

Taking into account  that the conditional distribution of the linear
form $L_2$ given $L_1$ is symmetric if and only if the conditional
distribution of the linear form $\alpha^{-1} L_2$
  given $L_1$ is symmetric,
  we may assume without loss of generality
  that $k\ge 0$, and hence the restriction of $\alpha \in {\rm
Aut}(\Omega_p)$ to the subgroup $\Delta_p$ is a continuous
endomorphism of  $\Delta_p$. We retain the notation $\alpha$ for
this restriction. We will construct in cases $2(i)$--$2(iv)$
independent random variables $\xi_1$ and $\xi_2$  with values in
$\Delta_p$ and distributions $\mu_1$ and $\mu_2$ such that the
conditional distribution of the linear form $L_2=\xi_1 + \alpha
\xi_2$ given $L_1=\xi_1+\xi_2$ is symmetric, whereas $\mu_1, \mu_2
\notin I(\Delta_p)$. Considering $\xi_j$ as independent random
variables with values in $\Omega_p$,  we prove statements
$2(i)$--$2(iv)$. Taking into account Lemma 1 it suffices to
construct non-idempotent distributions $\mu_j$ on the group
$\Delta_p$ such that their characteristic functions satisfy equation
(\ref{3n}). We note that the groups   ${\mathbf Z}(p^\infty)$ and
$\Delta_p$ are the character groups of one another, and the value of
a character $y={l\over p^n} \in {\mathbf Z}(p^\infty)$ at an element
$x=(x_0, x_1, \dots,x_n,\dots) \in \Delta_p$ is defined by the
formula
$$({x},
y)=\exp\displaystyle{\left\{\Bigl(x_0+x_1 p + \dots
+x_{n-1}p^{n-1}\Bigr){2\pi i l\over p^n}\right\}}.$$ Moreover, any
topological automorphism of the group $\Delta_p$ is the
multiplication by an element of  $\Delta^0_p.$ For $c=(c_0,
c_1,\dots, c_n,\dots)
 \in \Delta_p^0,$   the restriction of the automorphism $\widetilde c\in
{\rm Aut}({\mathbf{Z}}(p^\infty))$ to the subgroup
${\mathbf{Z}}(p^n) \subset {\mathbf{Z}}(p^\infty)$ is of the form
$\widetilde c y = s_n y,$ $y \in {\mathbf{Z}}(p^n)$, where
$s_n=c_0+c_1p+c_2p^2+\dots+c_{n-1}p^{n-1}$. Observe also that
$\widetilde\alpha=p^k\widetilde c$, and
\begin{equation}\label{15a}
A({\mathbf{Z}}(p^\infty), {p^l\Delta_p})={\mathbf{Z}}(p^l).
\end{equation}

Put $f(y)=\widehat\mu_1(y)$, $g(y)=\widehat\mu_2(y).$ In these
notation equation (\ref{3n})
 takes the form
\begin{equation}\label{29}
f(u+v)g(u+\widetilde\alpha v)=f(u-v)g(u-\widetilde\alpha v), \quad
u, v \in {\mathbf{Z}}(p^\infty).
\end{equation}

{\bf Statement $2(i)$}.
 Consider on the group $\Delta_p$
the distribution $\mu=a m_{\Delta_p}+(1-a)m_{p\Delta_p},$ where
$0<a<1$. It follows from (\ref{1a}) and (\ref{15a}) that the
characteristic function $\widehat\mu(y)$ is of the form
\begin{equation}\label{29a}
\widehat\mu(y)=
\left\{%
\begin{array}{ll}
    1, & \hbox{$y =0$;} \\
    1-a, & \hbox{$y \in
{\mathbf{Z}}(p)$;} \\
    0, & \hbox{$y\notin
{\mathbf{Z}}(p)$.} \\
\end{array}%
\right.
\end{equation}
Let us check that the characteristic functions
$f(y)=g(y)=\widehat\mu(y)$ satisfy equation (\ref{29}). Consider  3
cases:

1. $u, v \in {\mathbf{Z}}(p)$. Since $c_0=p-1$, we have
$\widetilde\alpha y = - y$ for $y \in {\mathbf{Z}}(p)$, and the
restriction of equation  (\ref{29}) to the subgroup
${\mathbf{Z}}(p)$ takes the form
\begin{equation}\label{30}
f(u+v)  g(u-v)=  f(u-v)  g(u+v), \quad u, v \in {\mathbf{Z}}(p).
\end{equation}
Since $f(y)=g(y)$,  (\ref{30}) holds.

2. Either $u \in {\mathbf{Z}}(p),$ $v \notin {\mathbf{Z}}(p)$, or $u
\notin {\mathbf{Z}}(p),$ $v \in {\mathbf{Z}}(p)$. Then $u \pm v
\notin {\mathbf{Z}}(p).$ This implies that $f(u \pm v)=0$,  and
hence both sides of equation (\ref{29}) are equal to zero.

3. $u, v \notin {\mathbf{Z}}(p)$. If $u+v, u+\widetilde\alpha v \in
{\mathbf{Z}}(p)$, then
\begin{equation}\label{31}
(I-\widetilde\alpha)v \in {\mathbf{Z}}(p).
\end{equation}
Since $p>2$, $k=0$ and $c_0=p-1$, we have $I-\alpha\in {\rm
Aut}(\Delta_p)$, and hence  $I - \widetilde\alpha \in {\rm
Aut}({\mathbf{Z}}(p^\infty)).$ Then  (\ref{31}) implies that $v \in
{\mathbf{Z}}(p)$, contrary to the assumption. Thus, either $u+v
\notin {\mathbf{Z}}(p)$ or $u+\widetilde\alpha v \notin
{\mathbf{Z}}(p)$, and the left-hand side of equation (\ref{29}) is
equal to zero. Similarly we check that the right-hand side of
equation (\ref{29}) is also equal to zero. Thus, (\ref{29}) holds.

\medskip

{\bf Statement $2(ii)$}. Consider  on the group $\Delta_2$ the
distribution $\mu=a m_{\Delta_2}+(1-a)m_{2\Delta_2},$ where $0<a<1$.
Then the characteristic function $\widehat\mu(y)$ is represented by
formula
 (\ref{29a}) for $p=2$. Let us check that the characteristic functions $f(y)=g(y)=\widehat\mu(y)$
satisfy equation (\ref{29}). Since $k=0$ and $c_0=c_1=1$, the
restriction of the  automorphism $\widetilde\alpha\in {\rm
Aut}({\mathbf{Z}}(2^\infty))$ to the subgroup ${\mathbf{Z}}(2^n)
\subset {\mathbf{Z}}(2^\infty)$ is of the form $\widetilde\alpha y =
m y,$ $y \in {\mathbf{Z}}(2^n)$, where
$m=1+2+c_22^2+\dots+c_{n-1}2^{n-1}=4l-1$. Consider 3 cases: 1. $u, v
\in {\mathbf{Z}}(2)$; 2. either $u \in {\mathbf{Z}}(2),$ $v \notin
{\mathbf{Z}}(2)$   or $u \notin {\mathbf{Z}}(2),$ $v \in
{\mathbf{Z}}(2)$; 3. $u, v \notin {\mathbf{Z}}(2)$. In cases 1 and 2
the reasoning is the same as in case $2(i)$.

Consider case 3, i.e. assume that $u, v \notin {\mathbf{Z}}(2)$, and
prove first that if $u+v, u+\widetilde\alpha v \in {\mathbf{Z}}(2)$,
then $u-v, u-\widetilde\alpha v \in {\mathbf{Z}}(2)$, and (\ref{29})
holds. To prove this note that the inclusions  $u+v,
u+\widetilde\alpha v \in {\mathbf{Z}}(2)$ are possible only in the
following cases.

$(a)$ $u+v=0$ and $u+\widetilde\alpha v=0$. This implies that
$(\widetilde\alpha-I)v=(m-1)v=2(2l-1)v=0$. Hence, $v \in
{\mathbf{Z}}(2)$, but this contradicts to the assumption.

$(b)$ $u+v={1\over 2}$  and $u+\widetilde\alpha v={1\over 2}$. The
reasoning is the same as in case $(a)$.

$(c)$ $u+v=0$ and $u+\widetilde\alpha v={1\over 2}$. This implies
that $(\widetilde\alpha-I)v=(m-1)v=2(2l-1)v={1\over 2}$. It follows
from this that either $v={1\over 4}$,  $u ={3\over 4}$, or
$v={3\over 4}$,  $u ={1\over 4}$. In both cases
 $u-v={1\over 2},$ $u-\widetilde\alpha v=u-mv=0$, and hence  (\ref{29}) holds.

$(d)$ $u+v={1\over 2}$, $u+\widetilde\alpha v=0$. The reasoning is
the same as in case $(c)$.

Reasoning similarly we verify that if  $u-v, u-\widetilde\alpha v
\in {\mathbf{Z}}(2)$, then $u+v, u+\widetilde\alpha v \in
{\mathbf{Z}}(2)$,  and (\ref{29}) holds.

\medskip

{\bf Statement $2(iii)$}.
  Consider  on the group $\Delta_2$
the distributions $\mu_1=a m_{2\Delta_2}+(1-a)m_{4\Delta_2}$ and
$\mu_2=a m_{\Delta_2}+(1-a)m_{2\Delta_2},$ where $0<a<1$. It follows
from (\ref{1a}) and (\ref{15a}) that that the characteristic
function $\widehat\mu_1(y)$ is of the form
$$
\widehat\mu_1(y)=
\left\{%
\begin{array}{ll}
    1, & \hbox{$y \in {\mathbf{Z}}(2)$;} \\
    1-a, & \hbox{$y \in {\mathbf{Z}}(4)$;} \\
    0, & \hbox{$y \notin {\mathbf{Z}}(4)$.} \\
\end{array}%
\right.
$$
Moreover,   the characteristic function $\widehat\mu_2(y)$ is
represented by formula
 (\ref{29a}) for $p=2$.  Let us check that the characteristic functions $f(y)=\widehat\mu_1(y)$
and $g(y)=\widehat\mu_2(y)$ satisfy equation (\ref{29}). Consider  3
cases:

1. $u, v \in {\mathbf{Z}}(4)$. Obviously, we can assume that
 $u \ne 0, \ v \ne 0$. Since either $c_0=1$, $c_1=0$ or $c_0=c_1=1$, the restriction of the automorphism
$\widetilde c$ to the subgroup ${\mathbf{Z}}(4)$ is of the form:
either $\widetilde y=y$  or $\widetilde y=-y$,  $y \in
{\mathbf{Z}}(4).$ Thus, the restriction of equation  (\ref{29}) to
the subgroup ${\mathbf{Z}}(4)$ either takes  the form
\begin{equation}\label{32}
  f(u+v)   g(u+2v)=   f(u-v)   g(u-2v), \quad
u, v \in {\mathbf{Z}}(4),
\end{equation}
or
\begin{equation}\label{33}
  f(u+v)  g(u-2v)=   f(u-v)   g(u+2v), \quad
u, v \in {\mathbf{Z}}(4).
\end{equation}
Consider equation  (\ref{32}). Equation (\ref{33})  can be
considered analogously.

$(a)$ $u=v={1\over 2}$. Then $u \pm v=0,$ $2v=0$. Hence, $ f(u \pm
v)=1,$ $ g(u \pm 2v)= g(u)$ and (\ref{32}) holds.

$(b)$ $u={1\over 2}$, $v \in \{{1\over 4}, {3\over 4}\}$. Then $ f(u
\pm v)= f(v)$. Since $u \pm 2v=0$, we have $  g(u \pm 2v)=1$, and
(\ref{32}) holds.

$(c)$ $u \in \{{1\over 4}, {3\over 4}\}$, $v ={1\over 2}$. Then $f(u
\pm v)=  f(u).$ Since $2v=0$, we have $  g(u \pm 2v)=  g(u)$, and
(\ref{32}) holds.

$(d)$ $u, v \in \{{1\over 4}, {3\over 4}\}$. This implies that $u
\pm 2v \notin {\mathbf{Z}}(2)$. Hence $  g(u \pm 2v)=0$, and  both
sides of equation (\ref{32}) are equal to zero.

2. Either $u \in {\mathbf{Z}}(4),$ $v \notin {\mathbf{Z}}(4)$ or $u
\notin {\mathbf{Z}}(4)$, $v \in {\mathbf{Z}}(4)$.  Then $u \pm
v\notin {\mathbf{Z}}(4).$ Hence $  f(u \pm v)=0$ and both sides of
equation (\ref{32}) are equal to zero.

3. $u, v \notin {\mathbf{Z}}(4)$. If  $u+v \in {\mathbf{Z}}(4)$ and
$u+\widetilde\alpha v \in {\mathbf{Z}}(2)$, then
\begin{equation}\label{16a}
  (I-\widetilde\alpha )v
\in {\mathbf{Z}}(4).
\end{equation}
 Since $k=1$, we have $I-\alpha\in {\rm Aut}(\Delta_2)$, and hence $I-\widetilde\alpha \in {\rm
Aut}({\mathbf{Z}}(2^\infty))$. Then (\ref{16a}) implies that $v \in
{\mathbf{Z}}(4)$, contrary to the assumption. Hence, either $u+v
\notin {\mathbf{Z}}(4)$ or $u+\widetilde\alpha v \notin
{\mathbf{Z}}(2).$ This implies that the left-hand side of equation
(\ref{29}) is equal to zero. Reasoning analogously, we verify that
the right-hand side of equation (\ref{29}) is also equal to zero.
Thus, (\ref{29}) holds.

\medskip

{\bf Statement $2(iv)$}.   Consider  on the group $\Delta_p$ the
distributions $\mu_1=a m_{p^{k-1}\Delta_p}+(1-a)m_{p^k\Delta_p}$ and
$\mu_2=a m_{\Delta_p}+(1-a)m_{p^{k-1}\Delta_p},$ where $0<a<1$. It
follows from (\ref{1a}) and (\ref{15a}) that   the characteristic
functions $\widehat\mu_1(y)$ and $\widehat\mu_2(y)$ are of the form
\begin{equation}\label{20b}
 \widehat\mu_1(y)=
    \left\{%
\begin{array}{ll}
    1, & \hbox{$y \in {\mathbf{Z}}(p^{k-1})$;} \\
    1-a, & \hbox{$y \in
{\mathbf{Z}}(p^k)\backslash{\mathbf{Z}}(p^{k-1})$;} \\
    0, & \hbox{$y \notin {\mathbf{Z}}(p^k)$.} \\
\end{array}%
\right.
    \quad\quad
    \widehat\mu_2(y) =
\left\{%
\begin{array}{ll}
    1, & \hbox{$y =
0$;} \\
    1-a, & \hbox{$y \in {\mathbf{Z}}(p^{k-1})\backslash\{0\}$;} \\
    0, & \hbox{$y \notin {\mathbf{Z}}(p^{k-1})$.} \\
\end{array}%
\right.
\end{equation}
Let us check that the characteristic functions
$f(y)=\widehat\mu_1(y)$ and $g(y)=\widehat\mu_2(y)$ satisfy equation
(\ref{29}). Consider  3 cases:

1. $u,  v \in {\mathbf{Z}}(p^k)$. Since $\alpha=p^kc$, we have
 $\widetilde\alpha v =
0$, and hence the restriction of equation (\ref{29}) to the subgroup
  ${\mathbf{Z}}(p^k)$ takes the form
\begin{equation}\label{17a}
 f(u+v)  g(u)=  f(u-v)   g(u), \quad u, v
\in {\mathbf{Z}}(p^k).
\end{equation}

$(a) \ u  \in {\mathbf{Z}}(p^{k-1})$. Then $  f(u \pm v)= f(v)$, and
equation (\ref{17a}) holds.

$(b) \ u \in  {\mathbf{Z}}(p^k)\backslash{\mathbf{Z}}(p^{k-1}).$
Then $  g(u)=0$, and both sides of equation (\ref{17a}) are equal to
zero.

2. Either $u \in {\mathbf{Z}}(p^k),$ $v \notin {\mathbf{Z}}(p^k)$ or
$u \notin {\mathbf{Z}}(p^k),$ $v \in {\mathbf{Z}}(p^k).$ Then  $u
\pm v \notin {\mathbf{Z}}(p^k)$. This implies that $f(u \pm v)=0$,
and both sides of equation (\ref{29}) are equal to zero.

 3.  $u, v \notin {\mathbf{Z}}(p^k).$
If $u+v \in {\mathbf{Z}}(p^k),$ $u+\widetilde\alpha v \in
{\mathbf{Z}}(p^{k-1})$, then
\begin{equation}\label{18a}
  (I-\widetilde\alpha )v
\in {\mathbf{Z}}(p^k).
\end{equation}
 Since $k\ge 1$, we have $I-\alpha\in {\rm Aut}(\Delta_p)$, and hence $I-\widetilde\alpha \in {\rm
Aut}({\mathbf{Z}}(p^\infty))$. Then (\ref{18a}) implies that $v \in
{\mathbf{Z}}(p^k)$, contrary to the assumption. Hence, either $u+v
\notin {\mathbf{Z}}(p^k)$ or $u+\widetilde\alpha v \notin
{\mathbf{Z}}(p^{k-1}).$ This implies that the left-hand side of
equation (\ref{29}) is equal to zero. Reasoning analogously, we
verify that the right-hand side of equation (\ref{29}) is also equal
to zero. Thus, (\ref{29}) holds. Theorem 1 is proved.

\medskip

Let $X=\Delta_p$. We remind that each automorphism $\alpha \in {\rm
Aut}(\Delta_p)$ is the multiplication by an element $c_\alpha \in
\Delta_p^0.$ Let  $\xi_1$ and $\xi_2$ be independent random
variables with values  in the group  $\Delta_p$. Then we can
consider $\xi_j$ as independent random variables with values  in the
group  $\Omega_p$. Moreover, it is obvious that the multiplication
by an element $c  \in \Delta^0_p$ is a topological isomorphism of
the group $\Omega_p$. This implies that statements $1(i)$ and
$1(ii)$ in Theorem 1 are also valid  for the group $\Delta_p$.
Taking into account that in the proof of statements $2(i)$ and
$2(ii)$ in Theorem 1 the corresponding independent random variables
take values in the subgroup $\Delta_p\subset\Omega_p$, we conclude
that the following theorem holds.

\medskip

{\bf Theorem 2}. {\it Let $X=\Delta_p$. Let $\alpha=c=(c_0, c_1,\dots,c_n, \dots)
 \in \Delta_p^0$ be an arbitrary topological automorphism of the group $X$.
 Then the following statements hold.

{$1$.} Let $\xi_1$ and $\xi_2$ be independent random variables with
values in $X$ and distributions $\mu_1$ and $\mu_2$. Assume that the
conditional distribution of the linear form $L_2=\xi_1 + \alpha
\xi_2$ given $L_1=\xi_1+\xi_2$ is symmetric.  Then

$1(i)$ If  $p>2$ and $c_0 \ne p-1$, then $\mu_j=m_K*E_{x_j}$, where
$K$ is a compact subgroup of
 $X$ and $x_j\in X$. Moreover, if $c_0=1$, then $\mu_1$ and $\mu_2$
are degenerate distributions.

$1(ii)$ If $p=2$,   $c_0=1$ and $c_1=0$, then $\mu_1$ and $\mu_2$
are degenerate distributions.

{$2$.} If one of the following conditions holds:

$2(i)$ $p>2,$   $c_0=p-1;$

$2(ii)$ $p=2,$   $c_0=c_1=1$,

\noindent then there exist independent random variables $\xi_1$ and
$\xi_2$  with values in $X$ and distributions $\mu_1$ and $\mu_2$
such that the conditional distribution of the linear form $L_2=\xi_1
+ \alpha \xi_2$ given $L_1=\xi_1+\xi_2$ is symmetric whereas $\mu_1,
\mu_2 \notin I(X)$.}

\medskip

{\bf Remark 1}. Assume that in Theorem 1 $p>2$, $k=0$, $c_0 \ne 1$
and  $c_0 \ne p-1$. This implies in particular, that $I-\alpha\in
{\rm Aut}(\Delta_p)$, and hence
\begin{equation}\label{19a}
 I-\widetilde\alpha\in {\rm Aut}({\mathbf{Z}}(p^\infty)).
\end{equation}
Let $\xi_1$ and $\xi_2$ be independent identically distributed
random variables with values in $\Delta_p$ and distribution
$m_{\Delta_p}$. We check that the conditional distribution of the
linear form $L_2=\xi_1 + \alpha \xi_2$ given $L_1=\xi_1+\xi_2$ is
symmetric. By Lemma 1, it suffices to verify that the characteristic
functions $f(y)=g(y)=\widehat m_{\Delta_p}(y)$ satisfy equation
(\ref{29}). It follows from (\ref{1a}) that
$$
\widehat m_{\Delta_p}(y) =
\left\{%
\begin{array}{ll}
    1, & \hbox{$y =0$;} \\
    0, & \hbox{$y \ne 0$.} \\
\end{array}%
\right.
$$
Obviously, it suffices to check that equation (\ref{29}) holds when
$u\ne 0$, $v\ne 0$. Thus, assume that   $u\ne 0$, $v\ne 0$.  If
$u+v=0$ and $u+\widetilde\alpha v=0$, then
$(I-\widetilde\alpha)v=0$. Taking into account (\ref{19a}), this
implies that $v=0$, contrary to the assumption. Thus, either $u+v\ne
0$ or $u+\widetilde\alpha v\ne 0$, and hence the left-hand side of
equation (\ref{29}) is equal to zero. Reasoning analogously we show
that the right-hand side of equation (\ref{29}) is also equal to
zero. So, equation (\ref{29}) holds. This example shows that
statement $1(i)$ in Theorem  1 can not be strengthened to the
statement that  both $\mu_1$  and $\mu_2$ are degenerate
distributions.

Assume that conditions $1(iii)$ of Theorem 1 hold, i.e. $p>2$ and
$|k|=1$. We can assume without loss of generality that $k=1$. Let
$\xi_1$ and $\xi_2$ be independent   random variables with values in
the group $\Delta_p$ and distributions
$\mu_1=am_{\Delta_p}+(1-a)m_{p\Delta_p}$, where $0<a<1$, and  $\mu_2
= m_{\Delta_p}$. Then the characteristic functions
$f(y)=\widehat\mu_1(y)$ and $g(y)=\widehat\mu_2(y)$ are defined by
formulas  (\ref{20b}) for $k=1$. Obviously, the proof that the
characteristic functions $f(y)$ and $g(y)$ satisfy equation
(\ref{29}), given in the proof of statement $2(iv)$, remains true
for $k=1$ too. Thus, statement $1(iii)$ can not be strengthened to
the statement that both distributions $\mu_1, \mu_2\in I(X).$

\medskip

{\bf Remark 2}. Compare Theorem 1 with Heyde's characterization theorem for two independent random variables on the real line
$X={\mathbf{R}}$.
It is easy to see that this theorem can be formulated in the following
way. Let  $\xi_1$ and $\xi_2$ be independent random variables with values in the group
$X={\mathbf{R}}$
and distributions $\mu_1$ and $\mu_2$. Let
$\alpha\ne -I$. If the conditional distribution of the linear form
$L_2=\xi_1 + \alpha
\xi_2$ given $L_1=\xi_1+\xi_2$ is symmetric, then $\mu_j$ are Gaussian. Thus, on the real line
the only condition $\alpha\ne -I$, obviously, necessary, is sufficient for characterization of
the Gaussian distribution.

It follows from Theorem 1 that on the group of $p$-adic numbers
$X=\Omega_p$ the state of affairs is more complicated.
On the one hand, much more  severe constrains for $\alpha$ (conditions $1(i)$ for $p>2$ and $1(ii)$ for $p=2$)
are necessary and sufficient
for characterization of the idempotent distribution. On the other hand,
there exist $\alpha$  (satisfying conditions $1(iii)$) such that the symmetry
of the conditional distribution of the linear form
$L_2=\xi_1 + \alpha
\xi_2$ given $L_1=\xi_1+\xi_2$ implies that the only one of the distributions $\mu_j$ is idempotent.
This effect is absent on the real line.

\end{document}